\documentclass[12pt]{amsart}
\usepackage{palatino}
\usepackage{amscd,amssymb}
\usepackage{float}
\usepackage{graphicx}

\usepackage[cmtip,matrix,arrow]{xy}

\usepackage{graphicx}

\numberwithin{equation}{section}

\newtheorem {Theorem}                   {Theorem}

\newtheorem {RefTheorem}[equation]      {Theorem}

\theoremstyle{definition}

\newtheorem {Question}[equation]      {Question}

 \newtheorem{zremark}[equation] {Remark}
\newenvironment{Remark} {\par\footnotesize\zremark}{~\\}{\endzremark}
\newcommand{\pr} {\smallskip\noindent{\bf Proof\,\,}}

\newcommand     {\comment}[1]   {}
\newcommand     {\mute}[2] {}
\newcommand     {\printname}[1] {}

\newcommand{\labell}[1] {\label{#1}\printname{#1}}

\def    \to     {\longrightarrow}

\def    \C      {{\mathbb C}}
\def    \R      {{\mathbb R}}

\def    \Z      {{\mathbb Z}}

\def    \pr     {\operatorname{pr}}

\begin{document}

\title[Lattice points in polytope boundaries]{Lattice points in polytope boundaries and formal geometric quantization of singular Calabi Yau hypersurfaces in toric varieties}

\author{Jonathan Weitsman}
\thanks{Supported in part by a grant from the Simons Foundation (\# 579801)}
\address{Department of Mathematics, Northeastern University, Boston, MA 02115}
\email{j.weitsman@neu.edu}
\thanks{\today}

\begin{abstract}  We show that the number of lattice points in the boundary of a positive integer dilate  of a Delzant integral polytope is a polynomial in the dilation parameter, analogous to the Ehrhart polynomial giving the number of lattice points in a lattice polytope.  We give an explicit formula for this polynomial, analogous to the formula of Khovanskii-Pukhlikov for the Ehrhart polynomial.  These counting polynomials satisfy a lacunarity principle, the vanishing of alternate coefficients, quite unlike the Ehrhart polynomial.  We show that formal geometric quantization of singular Calabi Yau hypersurfaces in smooth toric varieties gives this polynomial, in analogy with the relation of the Khovanskii-Pukhlikov formula to the geometric quantization of toric varieties.   The Atiyah-Singer theorem for the index of the Dirac operator gives a moral argument for the lacunarity of the counting polynomial.  We conjecture that similar formulas should hold for arbitrary simple integral polytope boundaries.

\end{abstract}

\maketitle

\tableofcontents

\section{Introduction} \subsection{The lattice point counting polynomial $R_{\partial\Delta}$ for polytope boundaries}  Let $\Delta \subset \R^n$ be an integral polytope in $\R^n:$  That is, a convex polytope with vertices on the integral lattice $\Z^n \subset \R^n.$  The polytope $\Delta$ is {\em simple} if the edge vectors at each vertex form a basis for $\R^n;$ we call such a polytope {\em Delzant} if, in addition, the edge vectors at each vertex can be chosen to form a basis for $\Z^n.$  Let $k \in \Z_{+}$ be a positive integer, and consider the expression

$$R_{\partial \Delta} (k) = \#(k \partial \Delta ) \cap \Z^n.$$

The function $R_{\partial \Delta} (k)$ is an analog of the Ehrhart polynomial

$$P_\Delta(k) = \#(k \Delta ) \cap \Z^n.$$

The fact that $P_\Delta(k)$ is a polynomial for any $\Delta,$ along with an inclusion-exclusion argument, shows that $R_{\partial \Delta} (k)$ is a polynomial as well.  Ehrhart-Macdonald reciprocity (see \cite{br}, Theorem 4.1) then implies the following.
\begin{RefTheorem}\labell{enum}  Let $\Delta$ be any lattice polytope.  Then
\begin{enumerate}\item The function $R_{\partial \Delta} (k) = \sum_{j=0}^{n-1} a_j k^{j} $ is a polynomial in $k$ of degree $n-1.$
\item$ a_{n-1} = {\rm vol} (\partial \Delta)$
\item $a_{n - 2j} = 0 $ for all $j = 1, \dots, [n/2].$
\end{enumerate} \end{RefTheorem}

In Part (2) of Theorem \ref{enum}, the volume of $\partial \Delta$ is defined as the sum of the volumes of the facets of $\Delta,$ each normalized to the sublattice generated by the facet vertices.  Alternatively, this is the sum of the Euclidean volumes of the facets, each divided by the length of the primitive lattice vector normal to it (See e.g. \cite{robins}, Lemma 5.18). The lacunarity principle (Part (3) in Theorem \ref{enum}) has no analog for the Ehrhart polynomial of a lattice polytope; but a similar lacunarity principle is observed in the weighted lattice point counts considered in Cappell-Shaneson \cite{cs}; see also \cite{cmss,ksw}.

The purpose of this paper is to derive an explicit formula for the polynomial $R_{\partial \Delta}(k)$ in the case that $\Delta$ is Delzant, and to understand the relation of this polynomial to the formal geometric quantization of singular Calabi Yau hypersurfaces in toric varieties, in analogy to the relation between the polynomial $P_\Delta$ and the geometric quantization of the toric variety associated to $\Delta$ (\cite{k,kk,kp, guillemin}).

To describe our formula for $R_{\partial \Delta}(k)$, we present $\Delta$ as an intersection of half-spaces in $\R^n.$

Let $d$ be the number of facets of $\Delta$, and let $n_1,\dots,n_d \in \Z^n$ be the primitive lattice vectors normal to the facets $F_i$ of $\Delta.$  Let $\lambda_1,\dots, \lambda_d \in \R_{+}$ be positive real numbers.  For each $i = 1,\dots, d,$ let 

$$H_i(\lambda_i) = \{ x \in \R^n:  \langle x , n_i \rangle \leq \lambda_i\}.$$

Then the Delzant polytope $\Delta$ may be written as

$$\Delta = \bigcap_{i=1}^d H_i(\mu_i)$$

\noindent where the $\mu_i$ are appropriately chosen positive real numbers.  This gives rise to a family of polytopes 

\begin{equation}\labell{d1}\Delta(\lambda_1,\dots,\lambda_d)  = \bigcap_{i=1}^d H_i(\lambda_i)\end{equation}

\noindent with 

\begin{equation}\labell{d2}\Delta(\mu_1,\dots,\mu_d) = \Delta.\end{equation}

Given $x \in \R$, we denote by $\hat{A}(x)$ the function

$$\hat{A}(x) = \frac{x/2}{\sinh{x/2}}.$$

The function $\hat{A}$ is analytic in a neighborhood of the origin, and has a power series 

\begin{equation}\labell{hata} \hat{A}(x) = 1 + \sum_{i=1}^\infty c_i x^{2i}.\end{equation}

Likewise, the function $\frac{1}{\hat{A}(x)}$ is analytic in a neighborhood of the origin, and has a power series

\begin{equation}\labell{hatai} \frac{1}{\hat{A}(x)} =  \frac{\sinh{x/2}}{x/2}= 1 + \sum_{i=1}^\infty d_i x^{2i}.\end{equation}

Then

\begin{Theorem}\labell{main2}  Let $\Delta \subset \R^n$ be a Delzant polytope given by an intersection of half spaces as in (\ref{d1}) and (\ref{d2}). The polynomial $R_{\partial \Delta} (k)$ is given by

\begin{equation}\labell{cpe}R_{\partial \Delta} (k) = (\prod_{i=1}^d \hat{A}(\frac{\partial}{\partial \lambda_i}))\frac{1}{\hat{A}}(\sum_{i=1}^d \frac{\partial}{\partial \lambda_i}){\rm vol} ({\partial \Delta}(\lambda_1,\dots,\lambda_d))\vert_{\lambda_i = k \mu_i} .\end{equation}
\end{Theorem}

\begin{Remark}
  In (\ref{cpe}) the infinite order differential operators  $\hat{A}(\frac{\partial}{\partial \lambda_i})$  are defined using the power series (\ref{hata}):

 $$\hat{A}(\frac{\partial}{\partial \lambda_i}):= \sum_j c_j (\frac{\partial}{\partial \lambda_i})^{2j}$$
 
 \noindent and similarly, using (\ref{hatai}),
 
 $$\frac{1}{\hat{A}}(\sum_{i=1}^d \frac{\partial}{\partial \lambda_i})= \sum_j d_j (\sum_{i=1}^d \frac{\partial}{\partial \lambda_i})^{2j}.$$
 
\noindent The volume ${\rm vol} ({\partial \Delta}(\lambda_1,\dots,\lambda_d))$ is defined as the sum of the Euclidean volumes of the facets of $\Delta(\lambda_1,\dots,\lambda_d),$  each divided by the length of the primitive lattice vector normal to it. Since the volume ${\rm vol} ({\partial \Delta}(\lambda_1,\dots,\lambda_d))$ is the volume of a union of polytopes given by the facets of $\Delta(\lambda_1,\dots,\lambda_d)$, it is a polynomial in the $\lambda_i,$ so that the expression (\ref{cpe}) is well-defined.  \end{Remark} 
 
 \subsection{Relation to Geometric Quantization of singular Calabi Yau hypersurfaces in toric varieties}  Let $\Delta \subset \R^n$ be a Delzant polytope as above.  Corresponding to $\Delta$ there exists an integral symplectic manifold $(M,\omega)$ of dimension $2n,$ equipped with an effective Hamiltonian action of the torus $T^n = (S^1)^n.$  The image $\phi(M)$ of the corresponding moment map $\phi : M \to \R^n = Lie(T^n)^*$ is the convex polytope $\Delta$ (see Delzant \cite{delzant}, and Guillemin \cite{guillemin}).  (In fact, we obtain symplectic forms $\omega_{\lambda_1,\dots,\lambda_d}$ on $M$ for ${\lambda_1,\dots,\lambda_d}$ sufficiently close to $\mu_1,\dots,\mu_d,$ and corresponding moment maps whose images are the polytopes $\Delta({\lambda_1,\dots,\lambda_d}).$)  The manifold $M$ is equipped with a complex line bundle $L$ and connection $\nabla$ of curvature $\omega.$  Furthermore, $M$ is a Kahler manifold, the line bundle $L$ is a holomorphic line bundle, giving a smooth polarized toric variety, and we may define the geometric quantization of the quadruple $(M,k \omega,L^k,\nabla_k)$ \footnote{Here $\nabla_k$ denotes the connection on $L^k$ arising from the connection $\nabla$ on $L$.} (where $k$ is a positive integer) in the Kahler polarization by
 
 $$Q_k(M) = H^0(M,L^k).$$
 
 Alternatively, we may consider the Dolbeault operator $\bar{\partial}_{L^k}$ and its index;\footnote{In this paper, by the index of an elliptic operator we mean both the virtual vector space consisting of the direct difference of the kernel and cokernel of the operator, and the dimension of that virtual vector space.  The intended meaning should be clear from the context.} we have
 
 $$ Q_k(M) = {\rm ind}(\bar{\partial}_{L^k}).$$
 
Then "quantization commutes with reduction" (in this case of a Kahler manifold, due to \cite{gs}) shows that

\begin{equation}\labell{gq}{\rm dim~}Q_k(M) = \# (k\Delta \cap \Z^n).\end{equation}

There is a further equivariant version of this theorem:  The quantization $Q_k(M)$ carries a torus representation, given by

\begin{equation}\labell{gq1}Q_k(M) = \oplus_{\alpha \in k \Delta\subset Lie(T^n)^*} \C_{\alpha}\end{equation}

\noindent where the sum is over weights $\alpha$ of the torus $T^n,$ and where each $\C_{\alpha}$ is the irreducible representation of $T^n$ given by the weight $\alpha.$  

\begin{Remark}\labell{gtvkp}{\bf The geometry of toric varieties and the Khovanskii-Pukhlikov formula.} We review the Khovanskii-Pukhlikov formula and the relation to Geometric Quantization (See e.g. \cite{guillemin}).  Let

$$Td(x) = \frac{x}{1-e^{-x}}.$$

\noindent The function $Td(x)$ is analytic near the origin and has the power series expansion

$$Td(x) = 1 + \sum_{k=1}^\infty \frac{b_k}{k!}  x^k$$

\noindent where the coefficients $b_k$ are (up to signs) the Bernoulli numbers (See e.g. \cite{bourbaki}).

We recall the following facts.
\begin{enumerate}
\item The Riemann-Roch Theorem

$$ {\rm ind}(\bar{\partial}_{L^k}) = \int_M Td(TM) e^{[\omega]}.$$

\item "Quantization Commutes with Reduction"

$$ {\rm ind}(\bar{\partial}_{L^k})  = \# (k\Delta \cap \Z^n).$$

\item The Duistermaat-Heckman theorem

$${\rm vol} (\Delta(\lambda_1,\dots,\lambda_d)) = \int_M e^{\omega_{\lambda_1,\dots,\lambda_d}}$$

\noindent and 

$$[\omega_{\lambda_1,\dots,\lambda_d}] =\sum_{i=1}^d  \lambda_i c_1(L_i)$$

\noindent where $L_i$ is the line bundle corresponding to the divisor given by $\phi^{-1}(F_i),$ and where $F_i$ is the $i$-th facet of $\Delta.$
\item The stable equivalence 

$$TM  \simeq \oplus_{i=1}^d L_i.$$

\end{enumerate}

Combining (3) and (4),we have

$$\prod_{i=1}^d (Td(\frac{\partial}{\partial \lambda_i})) \int_M e^{\omega_{\lambda_1,\dots,\lambda_d}} = \prod_{i=1}^d (Td(c_1(L_i)))e^{\omega_{\lambda_1,\dots,\lambda_d}} = \int_M Td(TM) e^{\omega_{\lambda_1,\dots,\lambda_d}}$$

\noindent and thus, using (1), (2), and (3), we obtain the Khovanskii-Pukhlikov formula \cite{k,kk,kp}

\begin{equation}\labell{kpff}{\rm dim~}Q_k(M) = \prod_{i=1}^d (Td(\frac{\partial}{\partial \lambda_i}))  {\rm vol~}(\Delta(\lambda_1,\dots,\lambda_d))\vert_{\lambda_i = k \mu_i} \end{equation}

\noindent where the infinite order differential operator $\prod_{i=1}^d (Td(\frac{\partial}{\partial \lambda_i}))$ is defined using the power series expansion at the origin of the function $Td(x),$ and is applied to the polynomial $ {\rm vol~}(\Delta(\lambda_1,\dots,\lambda_d)),$ as in Theorem \ref{main2}.
\end{Remark}

Consider now the locus

$$Z = \phi^{-1}(\partial \Delta) \subset M.$$

Then $Z$ is a singular subvariety of $M.$   

\begin{Remark}\labell{zcy} {\bf $Z$ as a singular Calabi Yau hypersurface of $M.$} Morally, $TZ \simeq TM|_Z \ominus NZ,$ where $NZ = \otimes_{i=1}^d L_i|_Z$ (where the $L_i$ are the line bundles on $M$ corresponding to the facets of $\Delta,$ as in Remark \ref{gtvkp}) is morally a normal bundle to $Z$ in $M.$\footnote{In the sense that the inverse image $\phi^{-1}(F_i)$ is the divisor corresponding to $L_i,$ so that, morally, the union $Z$ of these inverse images is the divisor corresponding to $\otimes_{i=1}^d L_i.$ Of course $Z$ is singular so this is only a moral argument, not a mathematical proof.}   
Note that then, since $TM \simeq \oplus_{i=1}^d L_i,$ we have 

$$c_1(TZ) = c_1(\oplus_{i=1}^d L_i|_Z) - c_1(\otimes_{i=1}^d L_i|_Z) = 0,$$ 

\noindent which is the sense in which $Z$ has vanishing first Chern class and may be viewed morally as a singular Calabi Yau hypersurface in $M.$   Batyrev \cite{batyrev} shows that where $\Delta$ satisfies a certain condition, $Z$ may be deformed to a Calabi Yau subvariety of $M$ with mild singularities.\end{Remark}

In analogy with (\ref{gq}) and (\ref{gq1}), we define the {\em formal Geometric Quantization} (see also \cite{weitsman, paradan}) of $Z$ by

$$Q_k(Z) = \oplus_{\alpha \in k \partial\Delta\subset Lie(T^n)^*} \C_{\alpha}$$

\noindent where the sum is again over weights $\alpha$ of the torus $T^n,$ and where each $\C_{\alpha}$ is the irreducible representation of $T^n$ given by the weight $\alpha.$  Then

\begin{equation}\labell{gqz}{\rm dim~}Q_k(Z) = \# (k\partial\Delta \cap \Z^n).\end{equation}

Then Theorem \ref{main2} may be stated as the following analog of (\ref{kpff}):

\begin{Theorem}\labell{gqmt} 
   Let $\Delta \subset \R^n$ be a Delzant polytope given by an intersection of half spaces as in (\ref{d1}) and (\ref{d2}). Let $k$ be a positive integer and let $Q_k(Z)$ be the formal geometric quantization of the singular Calabi Yau hypersurface $Z$ in the corresponding toric variety $M.$  Then 

\begin{equation}\labell{cpe2}{\rm dim~}Q_k(Z) = (\prod_{i=1}^d \hat{A}(\frac{\partial}{\partial \lambda_i}))\frac{1}{\hat{A}}(\sum_{i=1}^d \frac{\partial}{\partial \lambda_i}){\rm vol} ({\partial \Delta}(\lambda_1,\dots,\lambda_d))\vert_{\lambda_i = k \mu_i} .\end{equation}
\end{Theorem}
Note that again the volume ${\rm vol} ({\partial \Delta}(\lambda_1,\dots,\lambda_d))$ is defined as the sum of the Euclidean volumes of the facets of ${\partial \Delta}(\lambda_1,\dots,\lambda_d),$ each divided by the norm of the corresponding primitive lattice normal vector.
\begin{Remark}\labell{thm3rmk}{\bf The geometry of Calabi Yau manifolds and Theorem \ref{gqmt}.}  Note the lacunarity principle given by Part (3) of Theorem \ref{enum} is expected for an analog of geometric quantization of smooth Calabi Yau manifolds.  Suppose we are given a smooth Calabi Yau manifold $X$ of real dimension $2n,$\footnote{Note that the dimension of the singular Calabi Yau hypersurface $Z$ considered in this paper is taken to be $2n-2,$ two less than the dimension of the toric variety $M.$} along with a holomorphic line bundle $L.$  Then the Riemann Roch Theorem gives, for each integer $k,$

$${\rm ind~}(\bar{\partial}_{L^k}) = \int_X Td(TX) e^{k c_1(L)}$$

\noindent which is a polynomial in $k$ analogous to our $R_{\partial\Delta}(k).$

But since $X$ is Calabi Yau, it has trivial canonical class, and hence

$$Td(TX) = \hat{A}(TX) e^{\frac12 c_1(TX) } = \hat{A}(TX).$$

\noindent (See e.g. \cite{shan}, page 57.)

Thus 

\begin{equation}\labell{cyif} {\rm ind~}(\bar{\partial}_{L^k}) = \int_X \hat{A}(TX) e^{k c_1(L)}.\end{equation}

But the expression $\hat{A}(TX)$ is a cohomology class of mixed degree containing only even degree terms.  Hence, the expression ${\rm ind~}(\bar{\partial}_{L^k})$ is a polynomial in $k$  of degree $n$ with leading term 

$$ k^n \int_X   c_1(L)^n$$

\noindent and with all terms of degree $n-(2j -1)$ vanishing for $j = 1, \dots, [n/2].$

In a more speculative vein:  In the case of the singular Calabi Yau hypersurface $Z \subset M,$ as we saw in Remark \ref{zcy}, morally, $TZ \simeq TM|_Z \ominus NZ,$ where $NZ = \otimes_{i=1}^d L_i|_Z$ (where the $L_i$ are the line bundles on $M$ corresponding to the facets of $\Delta,$ as in Remark \ref{gtvkp}) is morally a normal bundle to $Z$ in $M.$\footnote{We reiterate that $Z$ is singular so this is only a moral argument, not a mathematical proof.}   
 Running an argument analogous to that of Remark \ref{gtvkp} for $Z$ instead of $M,$ we obtain an analog of equation (\ref{kpff}), with the ratio 
 $$\hat{A}(TZ) = \hat{A}(TM|_Z)/\hat{A}(NZ)= (\prod_{i=1}^d \hat{A}(L_i|_Z))/\hat{A}(\otimes_{i=1}^d L_i|_Z)$$
\noindent arising from (\ref{cyif}) giving rise to the operators on the right hand side of (\ref{cpe2}) replacing the Todd operators in the Khovanskii-Pukhlikov formula.  This gives a kind of explanation from geometric quantization, though not an alternative mathematical proof, for Theorem \ref{gqmt}.  It would be interesting to see if this argument could be worked up to a proof.  Perhaps the methods of Cappell-Shaneson \cite{cs},  which are applicable to singular spaces, would be useful in this context.  Another possible approach would be a version of "quantization commutes with reduction" for a torus action on a suitable class of symplectic stratified spaces; see e.g. \cite{sl}.

\end{Remark}

\subsection{Structure of the paper}  The main tool for proving Theorem \ref{main2} (as well as Theorem \ref{gqmt}) is the  the Khovanskii-Pukhlikov formula \cite{k,kk,kp} and its generalization by Agapito-Weitsman \cite{aw} (see also \cite{a,cmss,ms}) to weighted sums of lattice points.  These weighted sums include as a special case the number of interior points of the polytope $\Delta.$  In Section \ref{section2} we recall these formulas, and then subtract the number of interior points, as given by \cite{aw} from the number of all points of $\Delta,$ given by \cite{k,kk,kp}, to obtain a formula for the number of lattice points in $\partial \Delta.$  After some work we obtain formula (\ref{cpe}).  In Section \ref{section3} we first ask about the existence of a similar result for all simple integral polytopes (to which the formulas of \cite{cs,ksw,bv,guillemin} and \cite{aw} apply), corresponding to extending these formulas to toric varieties with orbifold singularities.   We also make some remarks about mirror symmetry.  Finally we perform some explicit calculations of the polynomial $R_{\partial \Delta}$ for the cases where $\Delta$ is a simplex in dimensions 2, 3, and 4.  The corresponding singular Calabi Yau varieties are a 2-torus with three pinched meridians, as a subset of $CP^2;$ a singular $K3$ surface in $CP^3;$ and a singular quintic in $CP^4.$

Acknowldgements:  I was prompted to write up this long-hibernating work by a question asked by V. Roubstov at the conference "Painleve equations, higher Teichmuller theory and quantization" held at the Polytechnic University of Catalonia in March 2025.  I would like to thank Prof. Roubstov, as well as Prof. M. Mazzocco and Prof. E. Miranda,  for a productive visit to Barcelona.  

I would also like to thank Matthias Beck, Ansgar Freyer, Matthias Schymura, and Sinai Robins for pointing out to me the proof of Theorem \ref{enum} using Ehrhart-Macdonald reciprocity, and for the reference to \cite{br}.  It is a pleasure to thank Prof. Robins particularly for a very helpful correspondence on the normalization of the volume of the boundary, and for many other helpful and encouraging comments which helped me improve the paper.

\section{The weighted Euler-Maclaurin formulas and the counting polynomial $R_{\partial\Delta}(k)$}\labell{section2}

In this section we prove Theorem \ref{main2} (and therefore also Theorem \ref{gqmt}).  As a consequence we also give an alternative proof, for the case of Delzant polytopes, of Theorem \ref{enum}.

We begin by recalling the theorems of Khovanskii-Pukhlikov \cite{k,kk,kp} and Agapito-Weitsman \cite{aw}.

Recall the description of the polytope $\Delta$ as an intersection of half-spaces.

Let $d$ be the number of facets of $\Delta$, and let $n_1,\dots,n_d\in \Z^n$ be the primitive outward normals to the facets $F_i$ of $\Delta.$  Let $\lambda_1,\dots, \lambda_d \in \R_{+}$ be positive real numbers.  For each $i = 1,\dots, d,$ let 

$$H_i(\lambda_i) = \{ x \in \R^n:  \langle x , n_i \rangle \leq \lambda_i\}.$$

Then the Delzant polytope $\Delta$ may be written as

$$\Delta = \bigcap_{i=1}^d H_i(\mu_i)$$

\noindent where the $\mu_i$ are appropriately chosen positive real numbers.
This gives rise to a family of polytopes 

\begin{equation}\labell{d3}\Delta(\lambda_1,\dots,\lambda_d)  = \bigcap_{i=1}^d H_i(\lambda_i)\end{equation}

\noindent with 

\begin{equation}\labell{d4}\Delta(\mu_1,\dots,\mu_d) = \Delta.\end{equation}

\begin{RefTheorem}\labell{k,kk,kp}[Khovanskii-Pukhlikov; Kantor-Khovanskii \cite{k,kk,kp}]  Let $k$ be a positive integer.  The number of lattice points in the Delzant polytope $k \Delta$ is given by

\begin{equation}\labell{kpf} \#(k \Delta \cap \Z^n) = \prod_{i=1}^d (Td(\frac{\partial}{\partial \lambda_i}))  {\rm vol~}(\Delta(\lambda_1,\dots,\lambda_d))\vert_{\lambda_i = k \mu_i} \end{equation}

\noindent where the infinite order differential operators $Td(\frac{\partial}{\partial \lambda_i})$ are given by substituting the first order differential operators $\frac{\partial}{\partial \lambda_i}$ into the power series at the origin of the function

$$Td(x) = \frac{x}{1-e^{-x}}$$

\noindent given by

$$Td(x) = 1 + \sum_{k=1}^\infty \frac{b_k}{k!}  x^k$$

\noindent where the coefficients $b_k$ are (up to signs) the Bernoulli numbers (See e.g. \cite{bourbaki}).

\end{RefTheorem}

Recall that since ${\rm vol~}(\Delta(\lambda_1,\dots,\lambda_d))$ is a polynomial in the $\lambda_i$, the expression on the right hand side of (\ref{kpf}) is well-defined.

In \cite{aw} Agapito and Weitsman proved a generalization of the formula of Khovanskii-Pukhlikov to weighted sums.  A special case of their result (in their notation, at $q=0$) is the following.

\begin{RefTheorem}\labell{aw} \cite{aw} Let $k$ be a positive integer.  The number of lattice points in the interior of the Delzant polytope $k \Delta$ is given by

\begin{equation}\labell{awf}\#(k (\Delta-\partial \Delta) \cap \Z^n) = \prod_{i=1}^d (Td(-\frac{\partial}{\partial \lambda_i})) {\rm vol~}(\Delta(\lambda_1,\dots,\lambda_d))\vert_{\lambda_i = k \mu_i} \end{equation}

\noindent where the infinite order differential operators $Td(-\frac{\partial}{\partial \lambda_i})$ are given by substituting the first order differential operators $\frac{\partial}{\partial \lambda_i}$ into the power series at the origin of the function

$$Td(-x) = \frac{-x}{1-e^{x}}.$$

\end{RefTheorem}

Again, since ${\rm vol~}(\Delta(\lambda_1,\dots,\lambda_d))$ is a polynomial in the $\lambda_i$, the expression on the right hand side of (\ref{awf}) is well-defined.

We therefore have, by subtraction

\begin{equation}\labell{subt}
\#(k \partial \Delta \cap \Z^n) = \left(\prod_{i=1}^d (Td(\frac{\partial}{\partial \lambda_i})) - \prod_{i=1}^d (Td(-\frac{\partial}{\partial \lambda_i})) \right) {\rm vol~}(\Delta(\lambda_1,\dots,\lambda_d))\vert_{\lambda_i = k \mu_i} .\end{equation}

To simplify the expression (\ref{subt}), we first note that

$$Td(-x) = e^{-x} Td(x).$$

Thus, for any $x_1,\dots,x_d \in (-1,1),$ we have

$$\prod_{i=1}^d Td(x_i) - \prod_{i=1}^d  Td(-x_i) = (1 - e^{-\sum_{i=1}^d  x_i}) \prod_{i=1}^d  Td(x_i).$$

But 

$$Td(x) = e^{x/2} \hat{A}(x)$$

so that

$$ (1 - e^{-\sum_{i=1}^d  x_i}) \prod_{i=1}^d  Td(x_i) = \prod_{i=1}^d  \hat{A}(x_i) (e^{\sum_{i=1}^d  x_i/2} - e^{-\sum_{i=1}^d  x_i/2})=
\prod_{i=1}^d  \hat{A}(x_i) \left(\frac{1}{\hat{A}(\sum_{i=1}^d  x_i)} \right) \sum_{i=1}^d  x_i.$$

Thus

$$\#(k \partial \Delta \cap \Z^n)=  (\prod_{i=1}^d \hat{A}(\frac{\partial}{\partial \lambda_i}))\frac{1}{\hat{A}}(\sum_{i=1}^d  \frac{\partial}{\partial \lambda_i})(\sum_{i=1}^d \frac{\partial}{\partial \lambda_i})  {\rm vol} ( \Delta(\lambda_1,\dots,\lambda_d))\vert_{\lambda_i = k \mu_i} .$$

Since 

$$\left(\sum_{i=1}^d \frac{\partial}{\partial \lambda_i}\right) {\rm vol} ( \Delta(\lambda_1,\dots,\lambda_d))=
{\rm vol} ( \partial \Delta(\lambda_1,\dots,\lambda_d)),$$

\noindent where ${\rm vol} ( \partial \Delta(\lambda_1,\dots,\lambda_d))$ is the sum of the volumes of the facets $F_i$ of $\Delta(\lambda_1,\dots,\lambda_d)),$ each normalized by the length $||n_i||$ of the primitive lattice normal $n_i,$\footnote{To see this, note that the variation of the Euclidean volume of a polytope as a facet is moved in the direction of a unit normal is given by the Euclidean volume of the facet.  Dilating the normal vector results in a corresponding change the the normalization of the volume of the facet.} we have

$$\#(k \partial \Delta \cap \Z^n)=  (\prod_{i=1}^d \hat{A}(\frac{\partial}{\partial \lambda_i}))\frac{1}{\hat{A}}(\sum_i \frac{\partial}{\partial \lambda_i})  {\rm vol} (\partial \Delta(\lambda_1,\dots,\lambda_d))\vert_{\lambda_i = k \mu_i} $$

\noindent proving Theorem \ref{main2}.

We also obtain an alternative proof of Theorem \ref{enum}.  To do so, note first that the function

$$ {\rm vol} (k \partial \Delta)$$ 

\noindent is a polynomial of degree $n-1$ in $k,$ so that all its derivatives are polynomials of degree less than or equal to $n-1.$  In particular the function $R_{\partial \Delta}(k)$ is a polynomial of degree $n-1.$ This proves Part (1) of Theorem \ref{enum}.  To study the terms in this polynomial, we note that the function

$$f(x_1,\dots,x_d)= \left( \prod_{i=1}^d \hat{A}(x_i)\right) \left(\frac{1}{\hat{A}(\sum_{i=1}^d x_i)} \right) $$

\noindent is an analytic function of the variables $x_1,\dots,x_d$ in a neighborhood of the origin in $\R^d,$ with

$$f(0,\dots,0) = 1$$

\noindent and with only even terms in its power series, given by a product of series arising from (\ref{hata}) and (\ref{hatai}).

The same is then true of the product of power series giving the differential operator entering into the right hand side of (\ref{cpe}).  

The fact that the lowest order term in this differential operator is the constant $1$ means that the leading term in the polynomial $R_{\partial\Delta}(k) $ is $ {\rm vol} (k \partial \Delta),$ proving Part (2) of Theorem \ref{enum}.

The fact that this differential operator contains only even order derivatives implies Part (3) of Theorem \ref{enum}.

\section{Some remarks and computations}\labell{section3}

\subsection{Simple polytopes and orbifold toric varieties}  So far we have considered the case where the polytope $\Delta$ is Delzant, corresponding to the case where the toric variety $M$ is smooth.  On the other hand the formulas of Khovanskii-Pukhlikov can be extended \cite{cs,bv,guillemin,dr, kk} to the case of simple integral polytopes, corresponding to the case where the toric variety $M$ is allowed to have orbifold singularities.  The results of \cite{aw} also apply in this case.  (The results of Cappell-Shaneson \cite{cs} apply to more singular toric varieties as well.)   

\begin{Question} What is the analog of Theorem \ref{main2} for simple lattice polytopes? 
\end{Question}

\subsection{Mirror Symmetry}  The paper of Batyrev \cite{batyrev} constructs, for each reflexive polytope $\Delta,$ a dual polytope $\Delta^*.$  The corresponding pair of Calabi Yau hypersurfaces is then shown by Batyrev to be a Mirror pair.

\begin{Question}\labell{mir}  Is there a relation between the formal geometric quantizations of the corresponding singular Calabi Yau hypersurfaces; equivalently, between the polynomials $R_{\partial \Delta}$ and $R_{\partial \Delta^*}?$
\end{Question}

In reference to Question \ref{mir}, we note that in general, mirror pairs do not come equipped with polarizing line bundles.  However, where the pair is presented as hypersurfaces in toric varieties, each presented as corresponding to a lattice polytope, the corresponding toric varieties, and therefore the hypersurfaces, are equipped with a line bundle.

We also note that the arithmetic genera, and in fact the elliptic genera, of mirror pairs are equal; see \cite{bl}.

\subsection{Some explicit calculations}  In this section we calculate the polynomial $R_{\partial \Delta}$ in the case where $\Delta$ is the unit lattice simplex in two, three, and four dimensions, corresponding to degenerations of a torus, a K3 surface, and a quintic.

\subsubsection{The degeneration of a torus in $CP^2$}

Take $M = CP^2.$  The moment image for the $T^2$ action on $CP^2$ is the unit right triangle $\Delta_2$ in $\R^2$ (See Figure 1).
\begin{figure} 
       \includegraphics[width=.5\textwidth]{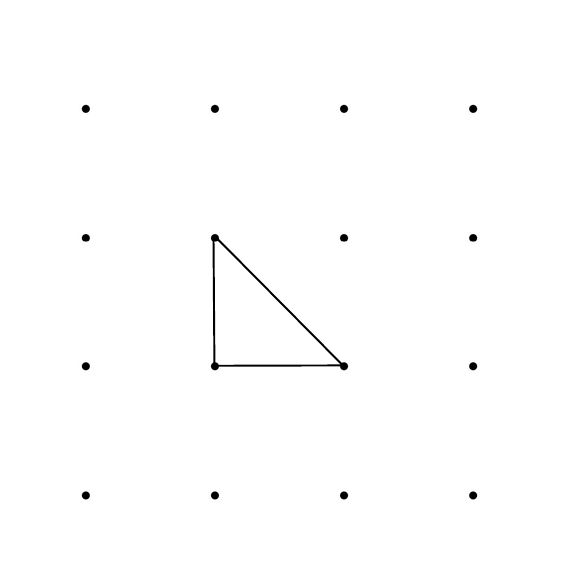}
\caption{The polytope $\Delta_2$}
   \end{figure}

The singular Calabi Yau hypersurface in this case is the singular curve $Z_1 = \{[x:y:z] \in CP^2 :  xyz = 0\}.$

By Theorem \ref{enum}, the polynomial $R_{\partial \Delta_2}(k)$ is given by

$$R_{\partial \Delta_2}(k) = 3k.$$

This can be verified by counting the lattice points on the boundary of the unit right triangle in $\Z^2,$ as in Figure 1.

Geometrically, the singular locus $Z_1$ in this case is the union of three copies of $CP^1$ glued to each other at three double points, corresponding to the degeneration of three meridians in the two torus.  See Figure 2.

\begin{figure} 
       \includegraphics[width=.5\textwidth]{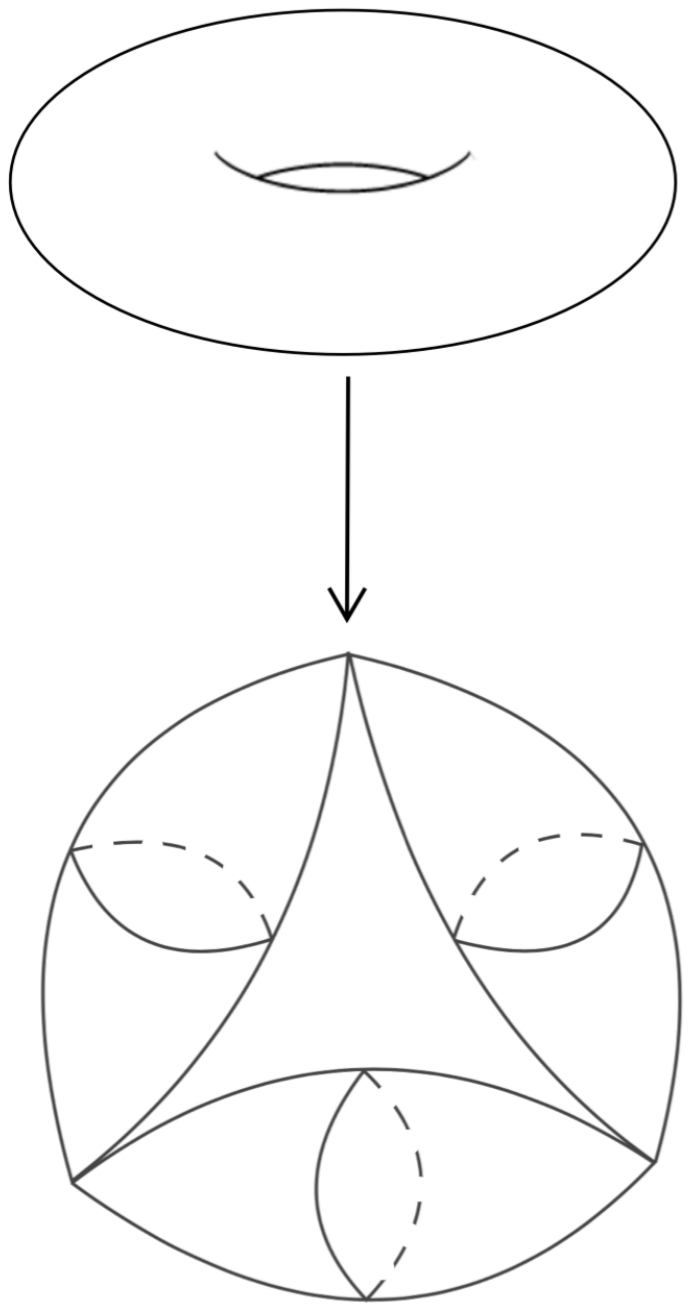}
\caption{The singular Calabi Yau $Z_1$ as a degeneration of a torus}
   \end{figure}

The absence of a constant term in the polynomial $R_{\partial \Delta_2}(k)$ is as expected from the Riemann Roch numbers of a torus. The integrality condition for the symplectic form on the degenerate locus $Z_1$ is responsible for the factor of $3$ in the expression for $R_{\partial \Delta_2}(k).$

\subsubsection{Singular K3 surface in $CP^3$}  Algebraically this is the locus

$$Z_2 = \{[x:y:z:t] \in CP^3 :  xyzt = 0\}.$$

\begin{figure} 
       \includegraphics[width=.5\textwidth]{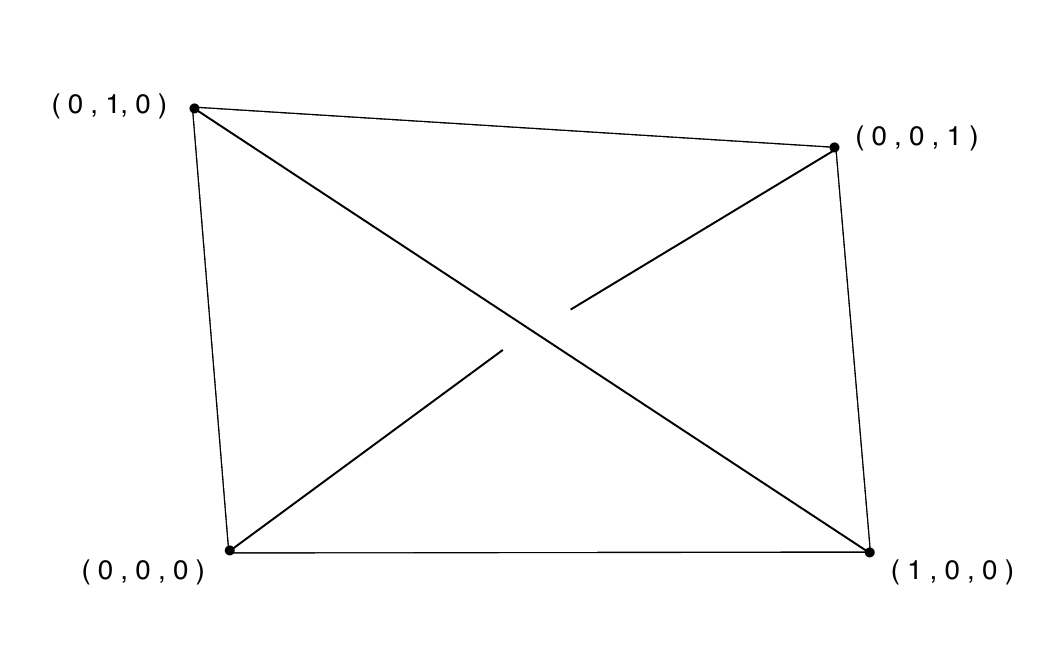}
\caption{The polytope $\Delta_3$}
   \end{figure}

The image of the moment map for the $T^3$ action on $CP^3$ is the unit tetrahedron $\Delta_3$  in $\R^3$ with
vertices $(0,0,0), (1,0,0), (0,1,0), (0,0,1).$  (See Figure 3.)

These vertices are the only lattice points on the boundary of the unit tetrahedron.  Thus the polynomial $R_{\partial \Delta_3}(k),$ which must have the form

$$R_{\partial \Delta_3}(k) = ak^2 + b$$

\noindent by the lacunarity principle, must have coefficients which satisfy

$$a + b = 4.$$

To determine $a$ and $b,$ we also count the number of lattice points in $2 \partial \Delta_3.$  The locus $2 \partial \Delta_3$ consists of four copies of the right triangle $2 \Delta_2.$  Thus there is one lattice point at each vertex of $\Delta_3,$ and the other lattice points are all located at the midpoint of each edge of $\Delta_3$  (See Figure 4.)  
\begin{figure} 
       \includegraphics[width=.5\textwidth]{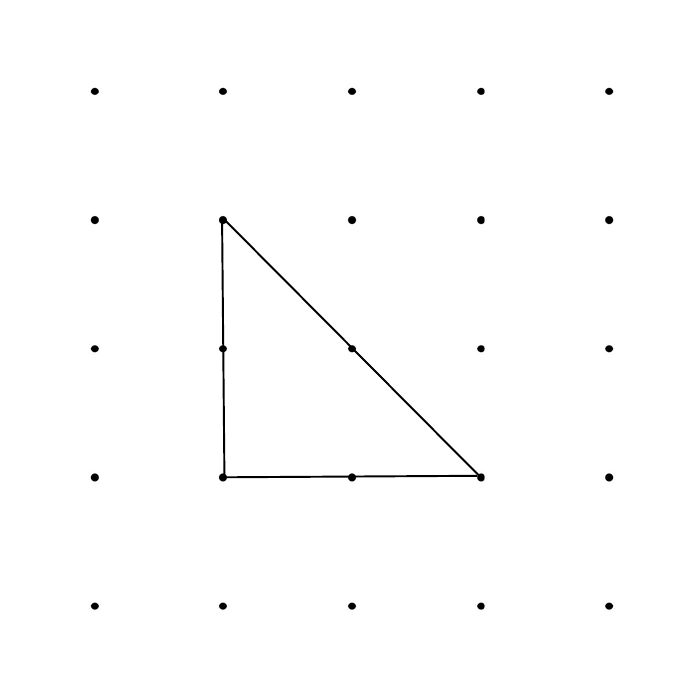}
\caption{The polytope $2\Delta_2$}
   \end{figure}

Hence

$$ 4 a + b = 10$$

\noindent so that 

$$R_{\partial \Delta_3}(k) = 2 k^2 + 2.$$

\begin{Remark}  It would be interesting to compare this lattice point count to the geometric quantization of $K3$ surfaces in \cite{mps}.\end{Remark}

\subsubsection{Singular quintic   in $CP^4$}  Algebraically this is the locus

$$Z_3 = \{[x:y:z:t:u] \in CP^4 :  xyztu = 0\}.$$

The image of the moment map for the $T^4$ action on $CP^4$ is the unit simplex $\Delta_4$  in $\R^4$ with
vertices $(0,0,0,0), (1,0,0,0), (0,1,0,0), (0,0,1,0),(0,0,0,1).$ These vertices are the only lattice points on the boundary of the unit tetrahedron.  (See Figure 5.)

\begin{figure} 
       \includegraphics[width=.5\textwidth]{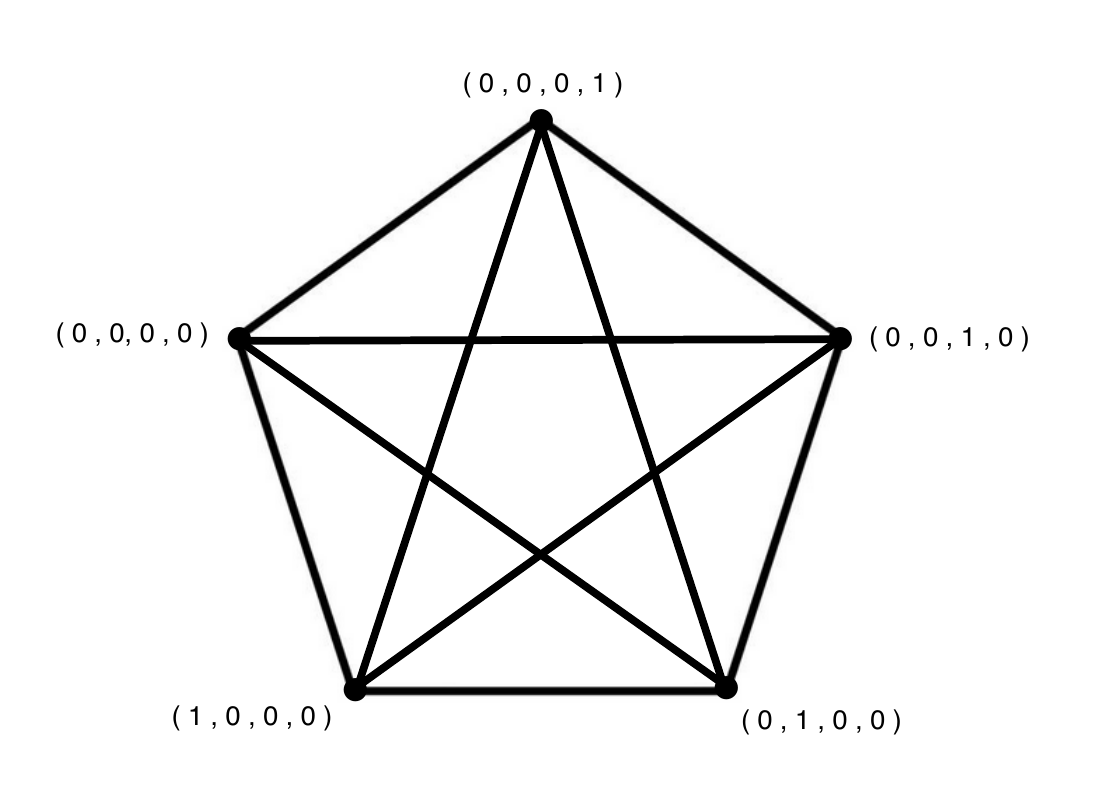}
\caption{The polytope $\Delta_4$}
   \end{figure}
Thus the polynomial $R_{\partial \Delta_4}(k),$ which must have the form

$$R_{\partial \Delta_4}(k) = ak^3 + bk$$

\noindent by the lacunarity principle, must have coefficients which satisfy

$$a + b = 5.$$

Again all the lattice points in $2\partial \Delta_4$ lie either at vertices or midpoints of edges, so that

$$ 8a + b = 15.$$

Thus 

$$R_{\partial \Delta_4}(k) = \frac56 k^3 + \frac{25}{6}k.$$

Note that  $k^3 + 5 k \equiv 0 ~~~({\rm mod~} 6)$
for all integers $k,$ so that $R_{\partial \Delta_4}(k)$ is an integer for each integer $k,$ as expected.

\end{document}